\documentclass[12pt,a4paper]{article}

\usepackage{anysize}
\usepackage[utf8]{inputenc}
\usepackage{amsmath}
\usepackage{amssymb}

\marginsize{3cm}{3cm}{3.5cm}{3.5cm}

\begin{document}
\title{On WZ-pairs which prove Ramanujan series}

\author{Jesús Guillera \\ jguillera@gmail.com \\ Av. Cesáreo Alierta, 31 esc. izda {\rm
$4^o$}--A, Zaragoza (Spain)}

\date{}

\maketitle

\newtheorem{observacion}{Remark}

\begin{abstract}
The known WZ-proofs for Ramanujan-type series related to $1/\pi$ gave us the insight to develop a new proof strategy based on the WZ-method. Using this approach we are able to find more generalizations and discover first WZ-proofs for certain series of this type.
\end{abstract}

Hypergeometric series. WZ-method 

\section{Introduction}
The Ramanujan series for $1/\pi$ are of the form \cite[pp. 352-354]{berndt}
\begin{equation}\label{rama-series}
\sum_{n=0}^{\infty} z^n \frac{ \left( \frac{1}{2} \right)_n (s)_n (1-s)_n}{ (1)_n^3} (a+bn)= \frac{1}{\pi},
\end{equation}
where $s$ determine the family according to the possible values $s=1/2$, $1/4$, $1/3$, $1/6$ and the parameters $z$, $a$, $b$ are algebraic real numbers with $-1 \leq z<1$. For an excellent survey on this kind of series, see \cite{baruah0}, which in addition cites many works of the main contributors. The symbol $(x)_n$ which appears in the series is the rising factorial or Pochhammer symbol, defined for all $x \in \mathbb{C}$ by
\begin{equation}\label{poch1}
(x)_n=\left\{
\begin{array}{ll}
x(x+1)\cdots(x+n-1), & \qquad n \in \mathbb{Z}^{+}, \\
1, & \qquad n=0,
\end{array} \right.
\end{equation}
or more generally by
\begin{equation}\label{poch2}
(x)_k=\frac{\Gamma(x+k)}{\Gamma(x)}.
\end{equation}
For $k \in \mathbb{Z}-\mathbb{Z}^{-}$, (\ref{poch2}) coincide with (\ref{poch1}). But (\ref{poch2}) is more general because it is also defined for all complex $x$ and $k$ such that
$x+k \in \mathbb{C}- ( \mathbb{Z}-\mathbb{Z}^{+} )$.
\par The Ramanujan series corresponding to rational values of the parameter $z$ can be written in the form
\begin{equation}\label{rama-racionales}
\sum_{n=0}^{\infty} z^n \frac{ \left( \frac{1}{2} \right)_n (s)_n (1-s)_n}{ (1)_n^3} (a+bn)= \frac{c}{\pi},
\end{equation}
where $a$ and $b$ are integers and $c^2$ is a rational number. D. Zeilberger in \cite{ekhad} proved the series for  $s=1/2$, $z=-1$, $a=1$, $b=4$, $c=2$ in a simple way, as an application of the WZ (Wilf and Zeilberger) method (Sect. \ref{wzmet}). Motivated by this result, in \cite{guilleraRJgen} and \cite{guilleratesis} we used WZ-pairs to obtain WZ-demonstrable generalizations of the series in Table I. In this paper we give a new strategy for the WZ-method specifically designed to prove Ramanujan-type series.
\begin{center}
\begin{tabular}{|c|c|c|c|c||c|c|c|c|c|}
  \hline
  $s$ & $z$ & $a$ & $b$ & $c$ & $s$ & $z$ & $a$ & $b$ & $c$  \\
  \hline \hline
  $1/2$ & $-1$     &  $1$   &  $4$   & $2$            &  $1/2$   & $1/4$      &   $1$    &  $6$    &  $4$           \\
  $1/2$ & $-1/8$   &  $1$   &  $6$   & $2 \sqrt{2}$   &  $1/2$   & $1/64$     &   $5$    &  $42$   &  $16$          \\
  $1/4$ & $-1/4$   &  $3$   &  $20$  & $8$            &  $1/4$   & $1/9$      &   $1$    &  $8$    &  $2/\sqrt{3}$  \\
  $1/4$ & $-1/48$  &  $3$   &  $28$  & $16/\sqrt{3}$  &  $1/6$   & $-27/512$  &   $15$   &  $154$  &  $32\sqrt{2}$  \\
  \hline
\end{tabular}
\end{center}
\begin{center} Table I \end{center}
Using this strategy, we have found other WZ-demonstrable generalizations of the series in Table I and the first WZ-proofs of those in Table II. We recall that the first proofs of the Ramanujan-type series for $1/\pi$ were based on the theory of modular forms and modular equations. For example, the formulas in Table I corresponding to $z=1/4$, $z=1/64$, $z=-1/4$, $z=-1/48$ and $z=1/9$ were discovered and proved in this way by S. Ramanujan \cite[Eq.: 28, 29, 35, 36, 40]{ramanujan}. See also
\cite[pp. 352-354]{berndt}.

\begin{center}
\begin{tabular}{|c|c|c|c|c|c|c|c|c|}
  \hline
  $s$ & $z$ & $a$ & $b$ & $c$ \\
  \hline \hline
  $1/3$  &  $1/2$      & $1$    & $6$     & $3\sqrt{3}$       \\
  $1/3$  &  $-9/16$    & $1$    & $5$     & $4/\sqrt{3}$      \\
  $1/3$  &  $-1/16$    & $7$    & $51$    & $12\sqrt{3}$      \\
  \hline
\end{tabular}
\begin{center} Table II \end{center}
\end{center}
The first proofs of the identities in Table II were also based on the theory of modular forms, but were discovered and proved much later \cite{chan}.

\section{The WZ-method}\label{wzmet}
We recall that a function $A(n,k)$ is hypergeometric in its two variables if the quotients
\[
\frac{A(n+1,k)}{A(n,k)} \quad {\rm and} \quad \frac{A(n,k+1)}{A(n,k)}
\]
are rational functions in $n$ and $k$, respectively. Also, a pair of hypergeometric functions $F(n,k)$ and $G(n,k)$ is said to be a Wilf and Zeilberger (WZ) pair \cite[Chapt. 7]{petkovsek} if
\begin{equation}\label{pro-WZ-pair}
F(n+1,k)-F(n,k)=G(n,k+1)-G(n,k).
\end{equation}
In this case, H. S. Wilf and D. Zeilberger \cite{wilf} proved that there exists a rational function $C(n,k)$ such that
\begin{equation}\label{certificado}
G(n,k)=C(n,k)F(n,k).
\end{equation}
The rational function $C(n,k)$ is the, so-called, certificate of the pair $(F,G)$. To discover WZ-pairs, usually we use a Maple package written by D. Zeilberger, called EKHAD \cite[Appendix A]{petkovsek}. If it certifies a function, we have found a WZ-pair!. However in this paper we only look for WZ-pairs which follow a special pattern, and we do not need that package. If we sum (\ref{pro-WZ-pair}) over all $n \geq 0$, we get

\begin{equation}\label{sumas-wz}
\sum_{n=0}^{\infty} G(n,k) - \sum_{n=0}^{\infty} G(n,k+1) = -F(0,k) + \lim_{n \to \infty} F(n,k).
\end{equation}
If $\lim_{n \to \infty} F(n,k)=0$ and $F(0,k)=0$, we have the identity
\[ \sum_{n=0}^{\infty} G(n,k) = \sum_{n=0}^{\infty} G(n,k+1). \]
If, in addition, the function
\[ f(k)=\sum_{n=0}^{\infty} G(n,k)-C, \]
where $C$ is a constant, satisfies the hypothesis of Carlson's theorem (see \cite{bailey}, p. 39), then we obtain
\[ \sum_{n=0}^{\infty} G(n,k)=C, \quad \forall k \in \mathbb{C}. \]
Usually it is easy to determine the constant by choosing a particular value for $k$. In other cases, we can determine the constant by taking the limit as $k \to \infty$.

\section{A new strategy}
We give a strategy to prove Ramanujan-type series
\begin{equation}\label{ramaform}
\sum_{n=0}^{\infty} z^n B(n) (a+bn)=\frac{c}{\pi},
\end{equation}
with $z$, $a$ and $b$ rational.
\par Let $B(n,k)$ be a function such that $B(n,0)=B(n)$ and hypergeometric in its two variables. We simplify the quotients $B(n+1,k)/B(n,k)$ and $B(n,k+1)/B(n,k)$ and denote the resulting denominators by $P(n,k)$ and $Q(n,k)$. Then we write
\[
P(n,k)=P_r(n,k) P_{r'}(n,k), \qquad Q(n,k)=Q_{s}(n,k) Q_{s'}(n,k),
\]
where $P_r(n,k)$ and $Q_{s}(n,k)$ are the polynomials of greatest possible degrees $r$ and $s$ satisfying $P_r(-1,0) \neq 0$ and $Q_{s}(0,-1) \neq 0$, respectively.
Finally, we construct two rational functions $R(n,k)$ and $S(n,k)$ (of degree 1) in the following way:
\begin{equation}\label{rama-RS}
R(n,k)=\frac{(a+bn) P_r(n,0) + k U_r(n,k)}{P_r(n,k)}, \qquad S(n,k)=\frac{n V_s(n,k)}{Q_s(n,k)},
\end{equation}
where
\[
U_r(n,k)=\sum_{j=0}^{r} \sum_{i=0}^{j} d_{i,j-i} n^i k^{j-i}, \qquad
V_s(n,k)=\sum_{j=0}^{s} \sum_{i=0}^{j} e_{i,j-i} n^i k^{j-i}
\]
are polynomials of degrees $r$ and $s$ whose coefficients must be determined.
A careful observation of the WZ-pairs used to prove Ramanujan-type series in \cite{guilleraRJgen} and \cite{guilleratesis} shows that they are of the form
\begin{equation}\label{rama-FG}
G(n,k)=z^n y^k B(n,k) R(n,k), \qquad F(n,k)=z^n y^k B(n,k) S(n,k),
\end{equation}
where $R(n,k)$ and $S(n,k)$ are rational functions of the form (\ref{rama-RS}) that we can construct in the way explained. To discover more WZ-pairs of this type, we consider other examples of functions $B(n,k)$ and try to solve for $y$, $d_{ij}$, $e_{ij}$ from the equation
\[ G(n,k+1)-G(n,k)=F(n+1,k)-F(n,k), \]
that is, from
\begin{equation}\label{eq2-WZ-pair}
\frac{B(n,k+1)}{B(n,k)}R(n,k+1)y-R(n,k)=\frac{B(n+1,k)}{B(n,k)} S(n+1,k)z-S(n,k).
\end{equation}
Simplifying (\ref{eq2-WZ-pair}), we obtain an equation of the form $H(n,k)=0$, where $H$ is a polynomial in $n$ and $k$. If there are values of $y$, $d_{ij}$ and $e_{ij}$ such that all the coefficients of $H(n,k)$ are equal to zero, then by substituting them in $F(n,k)$ and $G(n,k)$ we get a WZ-pair. In addition, as $F(0,k)=0$ we can apply the theory of the preceding section to prove that
\begin{equation}\label{general-rama}
\sum_{n=0}^{\infty} z^n y^k B(n,k) R(n,k)=\frac{c}{\pi}.
\end{equation}
In that case, the function $B(n,k)$ is good enough to prove a Ramanujan series.

\section{Types of searches}
Let us define the functions (by symmetry, $j_2=j_3$ when $s = 1/3, 1/4, 1/6$)
\begin{equation}\label{bin0}
C(n,k)=\frac{ \left( \frac{1}{2}+j_1 k \right)_n (s+j_2 k)_n ((1-s)_n+j_3 k)_n}{ (1)_n (1+j_4 k)_n (1_n+j_5 k)}, \qquad D(k)=\frac{ (t)_k (1-t)_k}{ (1)_k^2},
\end{equation}
where $t$ can be any of the numbers $1/2, 1/3, 1/4, 1/6$ and the $j_i$ are suitable rational numbers.
The simplest functions $B(n,k)$ that we have tried in this paper are of the form
\begin{equation}\label{bin1}
B(n,k)=C(n,k) D(k).
\end{equation}
In \cite{guilleraRJgen}, we found some series involving the simpler function
\begin{equation}\label{Dk-RJgen}
D(k)=\frac{\left( \frac{1}{2} \right)_k}{ (1)_k},
\end{equation}
and there are other candidates for $D(k)$ that we have not yet tried, for example
\begin{equation}\label{Dk-other}
D(k)=\frac{\left( \frac{1}{4} \right)_k \left( \frac{3}{4} \right)_k \left( \frac{1}{3} \right)_k \left( \frac{2}{3} \right)_k}{ \left( \frac{1}{2} \right)_k^2 (1)_k^2}.
\end{equation}
When we could not find solutions of the form (\ref{bin1}), we tried the longer expression
\begin{equation}\label{bin2}
B(n,k)=C(n,k) \cdot \frac{ \left( \frac{1}{2}+j_6 k \right)_n }{ \left( \frac{1}{2}+j_7 k \right)_n }
\cdot D(k),
\end{equation}
Of course, we can take more factors and so there are many other possibilities. For example, we can consider
\begin{equation}\label{bin3}
B(n,k)=C(n,k) \cdot \frac{ \left( \frac{1}{4}+j_6 k \right)_n \left( \frac{3}{4}+j_6 k \right)_n }
{ \left( \frac{1}{4}+j_7 k \right)_n \left( \frac{3}{4}+j_7 k \right)_n } \cdot D(k).
\end{equation}
However, we have not made the corresponding searches.
\\ \par Our searches are not completely random, because there are functions that we can discard in advance. For example,  the function
\begin{equation}\label{non-good1}
B(n,k)=\frac{ \left( \frac{1}{2}-k \right)_n \left( \frac{1}{3} \right)_n \left( \frac{2}{3} \right)_n}{ (1)_n (1+k)_n (1+2k)_n} \cdot \frac{\left( \frac{1}{3} \right)_k \left( \frac{2}{3} \right)_k}{(1)_k^2}
\end{equation}
cannot be good because it has a simple pole at $k=-1/2$, and so (\ref{general-rama}) cannot be finite at this value. Another example is the function
\begin{equation}\label{non-good2}
B(n,k)=\frac{ \left( \frac{1}{2} \right)_n \left( \frac{1}{2}+k \right)_n \left( \frac{1}{2}-5k \right)_n}{ (1)_n (1+k)_n^2} \cdot \frac{\left( \frac{1}{2} \right)_k^2}{(1)_k^2},
\end{equation}
which cannot be good, because at $k=\frac{1}{10}$ all the summands of (\ref{general-rama}) are zero except those corresponding to $n=0$ and we can easily check that we do not get a sum of the form $c/\pi$, where $c^2$ is a rational.

\section{An example}
Using our strategy, we prove in detail the following Ramanujan series in Table II:
\begin{equation}\label{rama-ejemplo}
\sum_{n=0}^{\infty} \left( \frac{-1}{16} \right)^n \frac{\left(
\frac{1}{2} \right)_n \left( \frac{1}{3} \right)_n \left(
\frac{2}{3} \right)_n}{(1)_n^3} (51n+7)=\frac{12 \sqrt{3}}{\pi}.
\end{equation}
\textbf{\textit{Solution 1:}} We will prove that the following
function $B(n,k)$ is good:
\begin{equation}\label{good1}
B(n,k)=\frac{\left( \frac{1}{2}-k \right)_n \left( \frac{1}{2}+k
\right)_n \left( \frac{1}{3} \right)_n \left( \frac{2}{3}
\right)_n}{\left(\frac{1}{2}+\frac{k}{2} \right)_n
\left(1+\frac{k}{2} \right)_n (1+k)_n (1)_n} \cdot \frac{\left(
\frac{1}{3} \right)_k \left( \frac{2}{3} \right)_k}{(1)_k^2}.
\end{equation}
\textit{Proof:} We define the functions
\[
G(n,k)=\left(\frac{-1}{16}\right)^n y^k B(n,k) R(n,k), \qquad
F(n,k)=\left(\frac{-1}{16}\right)^n y^k B(n,k) S(n,k),
\]
and get the evaluations
\begin{equation}\label{cocientes-B-1}
\frac{B(n+1,k)}{B(n,k)}=\frac{(2n-2k+1)(2n+2k+1)(3n+1)(3n+2)}{9(2n+k+1)(2n+k+2)(n+k+1)(n+1)},
\end{equation}
\begin{equation}\label{cocientes-B-2}
\frac{B(n,k+1)}{B(n,k)}=-\frac{(2n+2k+1)(3k+1)(3k+2)}{9(2n-2k-1)(2n+k+1)(n+k+1)}.
\end{equation}
Applying our strategy, we write
\[ R(n,k)=\frac{(51n+7)(2n+1)+k(d_{10}n+d_{01}k+d_{00})}{2n+k+1} \]
and
\[ S(n,k)= \frac{n(e_{10}n+e_{01}k+e_{00})}{2n-2k-1}. \]
Substituting in (\ref{eq2-WZ-pair}) and simplifying, we arrive at an equation of the form $H(n,k)=0$, where $H$ is a polynomial of degree $5$ and so has $21$ coefficients. As the coefficients of $n^5$, $n^4$, $n^3$, $n^2$, $n^4k$, $nk^4$ and of the independent term are, respectively,
\[ -58752+612e_{10}, \]
\[ -125568+1512e_{10}+612e_{00}, \]
\[ -60192+1187e_{10}+1476e_{00}-6528y, \]
\[ 23328+278e_{10}+1151e_{00}-7424y-64yd_{10}, \]
\[ -29376+1152e_{10}-576d_{10}+612e_{01}, \]
\[ -288yd_{10}+288d_{10}+144e_{01}-288yd_{01}+576d_{01}, \]
\[ 2016-2e_{10}-224y-2e_{00}-32yd_{01}-32yd_{00},  \]
we immediately get $e_{10}=96$, $e_{00}=-32$, $y=1$, $d_{10}=90$, $e_{01}=-48$, $d_{01}=24$, $d_{00}=28$.
Substituting these values in the other coefficients, we see that all of them are zero and so the polynomial is identically zero. Thus, $B(n,k)$ is a good function and the solution is given by $y=1$,
\begin{equation}\label{R-soluc1}
R(n,k)=\frac{(51n+7)(2n+1)+k(90n+24k+28)}{2n+k+1},
\end{equation}
and
\begin{equation}\label{S-soluc1}
S(n,k)= \frac{16n(6n-3k-2)}{2n-2k-1}.
\end{equation}
Then the WZ-method leads to
\[ \sum_{n=0}^{\infty} G(n,k)=\sum_{n=0}^{\infty} G(n,k+1), \]
and we can apply Carlson's theorem to derive
\[ \sum_{n=0}^{\infty} G(n,k)=C, \quad \forall k \in \mathbb{C}, \]
where we determine the constant $C$ by taking $k=1/2$. In this way, we have
the formula
\begin{multline}
\sum_{n=0}^{\infty} \left( \frac{-1}{16} \right)^n \frac{\left(
\frac{1}{2}-k \right)_n \left( \frac{1}{2}+k \right)_n \left(
\frac{1}{3} \right)_n \left( \frac{2}{3}
\right)_n}{\left(\frac{1}{2}+\frac{k}{2} \right)_n
\left(1+\frac{k}{2} \right)_n (1+k)_n (1)_n} \\ \times \frac{(51n+7)(2n+1)+k(90n+24k+28)}{2n+k+1} =
\frac{12 \sqrt{3}}{\pi} \cdot \frac{(1)_k^2}{\left( \frac{1}{3} \right)_k \left( \frac{2}{3} \right)_k}.
\end{multline}
Taking $k=0$ we obtain (\ref{rama-ejemplo}).
\\ \\ \textbf{\textit{Solution 2:}} Another good function $B(n,k)$ is
\begin{equation}\label{B-good2}
B(n,k)=\frac{\left( \frac{1}{2} \right)_n \left( \frac{1}{2}+2k
\right)_n \left( \frac{1}{3}+k \right)_n \left( \frac{2}{3}+k
\right)_n}{\left(\frac{1}{2}+\frac{k}{2} \right)_n
\left(1+\frac{k}{2} \right)_n (1+k)_n (1)_n} \cdot \frac{ \left(
\frac{1}{4} \right)_k \left(\frac{3}{4} \right)_k}{(1)_k^2}.
\end{equation}
\textit{Proof:} We define the functions
\[
G(n,k)=\left(\frac{-1}{16}\right)^n y^k B(n,k) R(n,k), \qquad
F(n,k)=\left(\frac{-1}{16}\right)^n y^k B(n,k) S(n,k).
\]
Applying our strategy, we obtain $y=1$,
\begin{equation}\label{R-soluc2}
R(n,k)=\frac{(51n+7)(2n+1)+k(114n+36k+37)}{2n+k+1},
\end{equation}
and
\begin{equation}
S(n,k)=\frac{-9n(6n^2+30nk+13n-7k-3)}{(3k+1)(3k+2)},
\end{equation}\label{S-soluc2}
The WZ method leads to
\[ \sum_{n=0}^{\infty} G(n,k)=\sum_{n=0}^{\infty} G(n,k+1), \]
and we can apply the Carlson's theorem to derive
\[ \sum_{n=0}^{\infty} G(n,k)=C, \quad \forall k \in \mathbb{C}, \]
where we determine the constant $C$ by substituting $k=-1/3$. We can
also find the constant by observing that as a consequence of the
Weierstrass M-test \cite[p. 49]{whittaker}, the convergence of the series is uniform. Hence the following steps hold:
\[
\lim_{k \to \infty} \sum_{n=0}^{\infty} G(n,k)=\sum_{n=0}^{\infty}
\lim_{k \to \infty} G(n,k)=\frac{18 \sqrt{2}}{\pi} \sum_{n=0}^{\infty} \left(
\frac{-1}{8} \right)^n {2n \choose n}=\frac{12 \sqrt{3}}{\pi},
\]
where we have used the identity \cite[p. 386]{borweinagm}
\[ \sum_{n=0}^{\infty} z^n {2n \choose n}=\frac{1}{\sqrt{1-4z}}. \]
Thus we have proved the formula
\begin{multline}\label{serie-gen-soluc2}
\sum_{n=0}^{\infty} \left( \frac{-1}{16} \right)^n \frac{\left( \frac{1}{2} \right)_n \left( \frac{1}{2}+2k \right)_n \left( \frac{1}{3}+k \right)_n \left( \frac{2}{3}+k \right)_n}{\left(\frac{1}{2}+\frac{k}{2} \right)_n \left(1+\frac{k}{2} \right)_n (1+k)_n (1)_n} \\ \times \frac{(51n+7)(2n+1)+k(114n+36k+37)}{2n+k+1} = \frac{12 \sqrt{3}}{\pi} \cdot \frac{(1)_k^2}{\left( \frac{1}{4} \right)_k \left( \frac{3}{4} \right)_k}.
\end{multline}
Taking $k=0$, we obtain (\ref{rama-ejemplo}).

\section{More formulas}\label{section-formulas}
For each of the indicated Ramanujan series (see Tables I and II), we have chosen an example among the WZ-demonstrable generalizations that we have found using our method.
\\ \\
Case: $s=1/2$, $z=-1$, $a=1$, $b=4$, $c=2$.
\begin{equation}\label{morefor1}
\sum_{n=0}^{\infty} (-1)^n \frac{\left( \frac{1}{2}-k \right)_n^2 \left( \frac{1}{2} \right)_n}{(1+k)_n^2 (1)_n}(4n+1)=\frac{2}{\pi} \cdot \left( \frac{1}{4} \right)^k \cdot \frac{(1)_k^2}{\left( \frac{1}{4} \right)_k \left( \frac{3}{4} \right)_k}.
\end{equation}
Case: $s=1/2$, $z=1/4$, $a=1$, $b=6$, $c=4$.
\begin{multline}\label{morefor2}
\sum_{n=0}^{\infty} \left( \frac{1}{4} \right)^n \frac{\left( \frac{1}{2}-k \right)_n \left( \frac{1}{2}+k \right)_n \left( \frac{1}{2}+3k \right)_n }{(1+k)_n (1+2k)_n (1)_n}(6n+6k+1) \\ = \frac{4}{\pi} \cdot \left( \frac{16}{27} \right)^k \cdot \frac{(1)_k^2}{\left( \frac{1}{6} \right)_k \left( \frac{5}{6} \right)_k}.
\end{multline}
Case: $s=1/2$, $z=-1/8$, $a=1$, $b=6$, $c=2 \sqrt{2}$.
\begin{multline}\label{morefor3}
\sum_{n=0}^{\infty} \left( \frac{-1}{8} \right)^n \frac{\left( \frac{1}{2}-k \right)_n \left( \frac{1}{2}+k \right)_n \left( \frac{1}{2}+3k \right)_n }{(1+k)_n (1+2k)_n (1)_n}(6n+6k+1) \\ = \frac{2 \sqrt{2}}{\pi} \cdot \left( \frac{32}{27} \right)^k \cdot \frac{(1)_k^2}{\left( \frac{1}{6} \right)_k \left( \frac{5}{6} \right)_k}.
\end{multline}
Case: $s=1/4$, $z=-1/4$, $a=3$, $b=20$, $c=8$.
\begin{multline}\label{morefor4}
\sum_{n=0}^{\infty} \left( \frac{-1}{4} \right)^n \frac{\left( \frac{1}{2}+k \right)_n \left( \frac{1}{4}+\frac{3k}{2} \right)_n \left( \frac{3}{4}+\frac{3k}{2} \right)_n }{(1+k)_n^2 (1)_n}(20n+18k+3) \\ = \frac{8}{\pi} \cdot \left( \frac{16}{27} \right)^k \cdot \frac{(1)_k^2}{\left( \frac{1}{6} \right)_k \left( \frac{5}{6} \right)_k}.
\end{multline}
Case: $s=1/4$, $z=1/9$, $a=1$, $b=8$, $c=2 \sqrt{3}$.
\begin{equation}\label{morefor5}
\sum_{n=0}^{\infty} \left( \frac{1}{9} \right)^n \frac{\left( \frac{1}{2} \right)_n \left( \frac{1}{4}+\frac{3k}{2} \right)_n \left( \frac{3}{4}+\frac{3k}{2} \right)_n }{(1+k)_n^2 (1)_n} (8n+6k+1) = \frac{2 \sqrt{3}}{\pi} \cdot \frac{(1)_k^2}{\left( \frac{1}{6} \right)_k \left( \frac{5}{6} \right)_k}.
\end{equation}
Case: $s=1/3$, $z=1/2$, $a=1$, $b=6$, $c=3 \sqrt{3}$.
\begin{equation}\label{morefor6}
\sum_{n=0}^{\infty} \left( \frac{1}{2} \right)^n \frac{\left( \frac{1}{2}+k \right)_n \left( \frac{1}{3} \right)_n \left( \frac{2}{3} \right)_n }{(1+k)_n (1+2k)_n (1)_n} (6n+6k+1) = \frac{3 \sqrt{3}}{\pi} \cdot  \frac{(1)_k^2}{\left( \frac{1}{3} \right)_k \left( \frac{2}{3} \right)_k}.
\end{equation}
Case: $s=1/2$, $z=1/64$, $a=5$, $b=42$, $c=16$.
\begin{multline}\label{morefor7}
\sum_{n=0}^{\infty} \left( \frac{1}{64} \right)^n \frac{\left( \frac{1}{2}-k \right)_n \left( \frac{1}{2} \right)_n \left( \frac{1}{2}+k \right)_n \left( \frac{1}{2}+2k \right)_n}{\left(\frac{1}{2}+\frac{k}{2} \right)_n
\left(1+\frac{k}{2} \right)_n (1+k)_n (1)_n} \\ \times \frac{(42n+5)(2n+1)+k(84n+24k+26)}{2n+k+1} = \frac{16}{\pi} \cdot \frac{(1)_k^2}{\left( \frac{1}{4} \right)_k \left( \frac{3}{4} \right)_k}.
\end{multline}
Case: $s=1/3$, $z=-1/16$, $a=7$, $b=51$, $c=12 \sqrt{3}$.
\begin{multline}\label{morefor8}
\sum_{n=0}^{\infty} \left( \frac{-1}{16} \right)^n \frac{\left(
\frac{1}{2} \right)_n \left( \frac{1}{2}+2k \right)_n \left(
\frac{1}{3}+k \right)_n \left( \frac{2}{3}+k \right)_n}{\left(\frac{1}{2}+\frac{k}{2} \right)_n
\left(1+\frac{k}{2} \right)_n (1+k)_n (1)_n} \\ \times \frac{(51n+7)(2n+1)+k(114n+36k+37)}{2n+k+1} = \frac{12 \sqrt{3}}{\pi} \cdot \frac{(1)_k^2}{\left( \frac{1}{4} \right)_k \left( \frac{3}{4} \right)_k}.
\end{multline}
Case: $s=1/3$, $z=-9/16$, $a=1$, $b=5$, $c=4/\sqrt{3}$.
\begin{multline}\label{morefor9}
\sum_{n=0}^{\infty} \left( \frac{-9}{16} \right)^n \frac{\left( \frac{1}{2}-k \right)_n \left( \frac{1}{2}+3k \right)_n \left( \frac{1}{3}+k \right)_n \left( \frac{2}{3}+k \right)_n}{\left(\frac{1}{2} \right)_n (1)_n (1+k)_n (1+3k)_n} \\ \times \frac{(5n+1)(2n+1)+k(16n+6k+7)}{2n+1} = \frac{4 \sqrt{3}}{3\pi} \cdot 4^k \cdot \frac{(1)_k^2}{\left( \frac{1}{6} \right)_k \left( \frac{5}{6} \right)_k}.
\end{multline}
Case: $s=1/4$, $z=-1/48$, $a=3$, $b=28$, $c=16 \sqrt{3}$.
\begin{multline}\label{morefor10}
\sum_{n=0}^{\infty} \left( \frac{-1}{48} \right)^n \frac{\left(
\frac{1}{2}-k \right)_n \left( \frac{1}{2}+3k \right)_n \left(
\frac{1}{4}\right)_n \left( \frac{3}{4} \right)_n}{\left(\frac{1}{2} \right)_n
(1)_n (1+k)_n^2} \\ \times \frac{(28n+3)(2n+1)+k(40n+18)}{2n+1} = \frac{16 \sqrt{3}}{\pi} \cdot \frac{(1)_k^2}{\left( \frac{1}{6} \right)_k \left( \frac{5}{6} \right)_k}.
\end{multline}
Case: $s=1/6$, $z=-27/512$, $a=15$, $b=154$, $c=32 \sqrt{2}$.
\begin{multline}\label{morefor11}
\sum_{n=0}^{\infty} \left( \frac{-27}{512} \right)^n \frac{\left( \frac{1}{2}-k \right)_n \left( \frac{1}{2}+k \right)_n \left( \frac{1}{6}+k \right)_n \left( \frac{5}{6}+k \right)_n}{\left(\frac{1}{2}+\frac{k}{2} \right)_n \left(1+\frac{k}{2} \right)_n (1+k)_n (1)_n} \\ \times \frac{(154n+15)(2n+1)+k(352n+108k+108)}{2n+k+1} \\ =
\frac{32 \sqrt{2}}{\pi} \cdot \left( \frac{32}{27} \right)^k \cdot \frac{(1)_k^2}{\left( \frac{1}{6} \right)_k \left( \frac{5}{6} \right)_k}.
\end{multline}
We have checked that the WZ-pairs used in \cite{guilleraAAMwz}, \cite{guilleraRJgen} and \cite{guilleratesis} to prove similar series for $1/\pi^2$ satisfy the same criteria, except that the degree of the functions $R(n,k)$ and $S(n,k)$ is now $2$. So we can try to use this strategy to prove the other unproved series in \cite{guilleraEMpi2}.

\section{ADDENDUM: A Maple program (July, 2019)}
I have written a Maple program which implements a generalization of the above theory (you can see it in my web-site). To use the program, begin writing the degree (grado), which is equal to $j$ for the Ramanujan-type series for $1/\pi^j$. Then write \texttt{bb:=(n,k)->} followed by the corresponding expression of pochhammer's symbols written in the Maple language. Finally, write \texttt{resolver(z,a,b,c,d)}, replacing the parameters with their corresponding values (when the degree is $1$, put $c=d=0$ and when it is $2$, put $d=0$). Thus, for example for (\ref{good1}), write

\begin{verbatim}
grado:=1;
bb:=(n,k)->po(1/2-k,n)*po(1/2+k,n)*po(1/3,n)*po(2/3,n)*po(1/3,k)*
po(2/3,k)/(po(1/2+k/2,n)*po(1+k/2,n)*po(1+k,n)*po(1,n)*po(1,k)^2);
resolver(-1/16,7,51,0,0);
\end{verbatim}
You can check the examples of Sect. {\bf \ref{section-formulas}}. Below we give examples of degree $2$.
\begin{verbatim}
grado:=2;

bb:=(n,k)->po(1/2,n)^5/(po(1,n)*po(1+k,n)^4)*po(1/2,k)^4/po(1,k)^4;
resolver(-1/4,1,8,20,0);

bb:=(n,k)->po(1/2,n)*po(1/2+k,n)^4/po(1,n)^5;
resolver(-4,1,6,10,0);

bb:=(n,k)->po(1/2,n)^3*po(1/4-k/2,n)*po(3/4-k/2,n)/
(po(1,n)^3*po(1+k,n)^2)*po(1/2,k)^2/po(1,k)^2;
resolver(1/16,3,34,120,0);

bb:=(n,k)->po(1/2,n)*po(1+k/2,n)^4*po(1/2+k/2,n)^4/
(po(1,n)^5*po(1+k,n)^4);
resolver(-1024,32,160,205,0);

bb:=(n,k)->po(1/2+k,n)^5*po(1/2-k,n)^4*po(1/2,k)/
(po(1,n)^4*po(1/2,n)^4*po(1+k,n)*po(1,k));
resolver(-1/1024,13,180,820,0);

bb:=(n,k)->po(1/2-k,n)^3*po(1/2+k,n)^3*po(1/3,n)*po(2/3,n)/
(po(1,n)^5*po(1/2,n)^3);
resolver(27/64,3,27,74,0);
\end{verbatim}
You can also check the examples of other papers or try to discover new ones.


\begin{thebibliography}{99}

\bibitem{bailey} Bailey, W.N.: Generalized Hypergeometric Series. Cambridge Univ. Press, (1935).

\bibitem{baruah0} Baruah, N.D., Berndt B.C., Chan, H.H.: Ramanujan's series for $1/\pi$: A survey,
The Amer. Math. Monthly \textbf{116} (2009) 567-587.; also available at the pages http://www.math.ilstu.edu/cve/speakers/Berndt-CVE-Talk.pdf, and http://www.math.uiuc.edu/~berndt/articles/monthly567-587.pdf.

\bibitem{berndt} Berndt, B.C.: Ramanujan's Notebooks, Part IV. Springer-Verlag, New York, (1994).

\bibitem{borweinagm} Borwein, J.M., Borwein, P.B.: Pi and the AGM: A Study in Analytic Number Theory and
Computational Complexity, (Canadian Mathematical Society Series of Monographs and Advanced Texts), Jonh Wiley, New York, (1987).

\bibitem{chan} Chan, H.H., Liaw, W.C., Tan, V.: Ramanujan's class invariant $\lambda_n$ and a new class of series
for $1/\pi$. J. London Math. Soc. \textbf{64}, 93-106, (2001).

\bibitem{ekhad} Ekhad, S.B., Zeilberger D.: A WZ proof of Ramanujan's formula for $\pi$. In Rassias, J.M. (ed.).
Geometry, Analysis and Mechanics. World Scientific, Singapore, 107-108, (1994); also available at arXiv:math/9306213v1. (The coauthor EKHAD is a Maple package written by D. Zeilberger).

\bibitem{guilleraAAMwz} Guillera, J.: Some binomial series obtained by the WZ-method.
Adv. in Appl. Math. \textbf{29}, 599-603, (2002).

\bibitem{guilleraEMpi2} Guillera, J.: About a new kind of Ramanujan type series.
Exp. Math. \textbf{12}, 507-510, (2003).

\bibitem{guilleraRJgen} Guillera, J.: Generators of Some Ramanujan Formulas. Ramanujan J. \textbf{11}, 41-48, (2006).

\bibitem{guilleratesis} Guillera, J.: Series de Ramanujan: Generalizaciones y
conjeturas. Ph.D. Thesis, University of Zaragoza, Spain, (2007).

\bibitem{petkovsek} Petkov\u{s}ek, M., Wilf, H.S., Zeilberger, D.: A=B. Peters A.K.: Ltd., (1996);
also available at http://www.math.upenn.edu/~wilf/AeqB.html.

\bibitem{ramanujan} Ramanujan, S.: Modular equations and approximations to $\pi$.
Q. J. Math. \textbf{45}, 350-372, (1914).

\bibitem{wilf} Wilf, H.S., Zeilberger, D.: Rational functions certify combinatorial identities.
Journal Amer. Math. Soc. \textbf{3}, 147-158, (1990). (Winner of the Steele prize).

\bibitem{whittaker} Whittaker, E.T., Watson, G.N.: A Course of Modern Analysis. Cambridge Univ. Press, (1927).

\end{thebibliography}
\end{document}